\documentclass[11pt,oneside,a4paper,final]{article}
\usepackage{amsmath,amssymb,makeidx,amsfonts,latexsym,enumerate,amsthm}
\usepackage{cite}

\newtheorem{lemma}{Lemma}[section]
\newtheorem{theorem}{Theorem}[section]

\newtheorem{proposition}{Proposition}[section]

\newtheorem{remark}{Remark}[section]

\newtheorem{example}{Example}[section]

\begin{document}

\author{Karl-Olof Lindahl\\
School of Mathematics and Systems Engineering\\
V\"{a}xj\"{o} University, 351 95, V\"{a}xj\"{o}, Sweden\\
\texttt{karl-olof.lindahl@lnu.se}}

\title{Linearization in ultrametric dynamics in fields of characteristic zero -- equal characteristic case\footnote{\textbf{Published in \emph{p-Adic Numbers, Ultrametric Analysis and Applications}}, Vol. 1, No 4, pp. 307--316, 2009}}

\date{October 13, 2009}

\maketitle

\begin{abstract}
Let $K$ be a complete ultrametric field of charactersitic zero whose corresponding residue field $\Bbbk$ is also of charactersitic zero.  We give lower and upper bounds for the size of linearization disks for
power series over $K$ near an indifferent fixed point. These estimates are maximal in the sense that there exist exemples where these estimates give the exact size of the corresponding linearization disc. Similar estimates in the remaning cases, i.e. the cases in which $K$ is either a $p$-adic field or a field of prime characteristic, were obtained in various papers on the $p$-adic case  
\cite{Ben-Menahem:1988,Thiran/EtAL:1989,Pettigrew/Roberts/Vivaldi:2001,Khrennikov:2001a} later generalized in \cite{Lindahl:2004cpp}, and in \cite{Lindahl:2004,Lindahl:2009} concerning the prime characteristic case.

\end{abstract}

\vspace{1.5ex}\noindent {\bf Mathematics Subject Classification
(2000):} 32P05, 32H50, 37F50

\vspace{1.5ex}\noindent {\bf Key words:} dynamical system,
linearization, conjugacy, ultrametric field

\section{Introduction}

A central issue in the study of dynamical systems is the local dynamics near periodic points. In particular, it is of great importance to know wether or not the  dynamics is locally linearizable near a given periodic point. Recall that a power series
$f$, over a complete valued field $K$, of the form
\begin{equation}\label{eq form of f}
f(x)=\lambda x+a_2x^2+a_3x^3\dots, \quad \text{with } |\lambda |=1, \text{but not a root of unity}, 
\end{equation}
is said to be analytically linearizable at the indifferent fixed point  at the origin if there is a convergent power series 
solution $g$ to the following form of the Schr\"{o}der functional equation
\begin{equation}\label{schroder functional equation}
g\circ f(x)=\lambda g(x), \quad \lambda =f'(0),
\end{equation}
which conjugates $f$ to its linear part. 
The coefficients of the formal solution $g$ of (\ref{schroder
functional equation}) must satisfy a recurrence relation of the form
\begin{equation*}
b_k=\frac{1}{\lambda (1-\lambda^{k-1})}C_k(b_1,\dots,b_{k-1}).
\end{equation*}
Intuitively, if $\lambda $ is close to a root of unity we might run into a problem of small divisors as in the well-known complex case \cite{Carleson/Gamelin:1991,Beardon:1991,Milnor:2000}. In 1942 Siegel proved in his celebrated paper
\cite{Siegel:1942} that the condition
\begin{equation}\label{Siegel condition}
|1- \lambda^n|\geq Cn^{-\beta} \quad\text{for some real numbers
$C,\beta >0$},
\end{equation}
on $\lambda$ is sufficient for convergence in the complex field
case. Later,  Brjuno \cite{Brjuno:1971} proved that the weaker
condition
\begin{equation}\label{condition brjuno}
-\sum_{k=0}^{\infty} 2^{-k}\log\left (\inf_{1\leq n\leq 2^{k+1}-1}
|1 - \lambda ^{n}| \right) <+\infty,
\end{equation}
is sufficient. 
In fact, for quadratic polynomials, the Brjuno condition is not only sufficient
but also necessary as shown by Yoccoz \cite{Yoccoz:1988}.

Since then, there has been an increasing interest in the ultrametric analogue
of complex dynamics, see e.g.\@ \cite{Arrowsmith/Vivaldi:1993,Arrowsmith/Vivaldi:1994,Benedetto:2001b,Benedetto:2003a,
Bezivin:2004b,Hsia:2000,Khrennikov:2001a,KhrennikovNilsson:2001,
Li:1996a,Li:2002a,Lindahl:2004,Lindahl:2007,Lubin:1994,Rivera-Letelier:2003thesis,Rivera-Letelier:2003,DeSmedtKhrennikov:1997,
AnashinKhrennikov:2009,Khrennikov/Nilsson:2004,
Khrennikov:2003ryssbok,DragovichKhrennikovMihajlovic:2007,
KhrennikovMukhamedovMendes:2007,KhrennikovSvensson:2007,
Khrennikov:2003nauk, Svensson:2005,NilssonNyqvist:2004}.

It is known since a work of Herman and Yoccoz \cite{HermanYoccoz:1981}
that Siegel's linearization theorem is true also in the ultrametric case. Moreover, in the one-dimensional case, for fields of characteristic zero the Siegel condition is always satisfied. In the two-dimensional $p$-adic case, the conjugacy may diverge as shown in\cite{HermanYoccoz:1981}.  Recently, the multi-dimensional $p$-adic case has been studied in more detail by Viegue in his thesis \cite{Viegue:2007}.

However, as noted by Herman and Yoccoz, in fields of prime characteristic there is a problem of small divisors also in the one-dimensional case; in general the multiplier does not satisfy the Siegel nor the weaker Brjuno condition. 
 One might therefore conjecture,
as Herman \cite{Herman:1986}, that
\emph{for a locally compact, complete valued field of prime 
characteristics, the formal conjugacy `usually'  diverges, even for polynomials
of one variable}.
Indeed, as shown in the papers \cite{Lindahl:2004,Lindahl:2009} 
like in complex dynamics, the formal solution may diverge also in the one-dimensional case.
On the other hand, in \cite{Lindahl:2004,Lindahl:2009} it was also proven that the
conjugacy may still converge due to considerable cancellation of
small divisor terms; the same multipler $\lambda$ may yield convergence for some $f$ but not for others. This brings about a problem of a
combinatorial nature of seemingly great complexity and a complete
description is yet to be found. For example, we have the following open problem stated in \cite{Lindahl:2004}. 

\vspace{1.5ex}
\noindent
\textbf{Open Problem}
\textit{Let $K$ be of characteristic $p>0$.
Is there a polynomial of the form $f(x)=\lambda x+
O(x^2)\in K[x]$, with $\lambda $ not a root of
unity satisfying $|1-\lambda ^n|<1$ for some $n\geq 1$, and
containing a monomial of degree prime to $p$, such that the formal
conjugacy $g$ converges?}\\

For a more thorough treatment of the problem and its relation to the complex case, the reader can consult \cite{Lindahl:2004}.

In case of convergence, one can estimate the radius of 
convergence for
the corresponding
 \emph{linearization disc}\footnote{Here we use the term `linearization disc' rather than `Siegel disc', because in ultrametric dynamics the Siegel disc is 
often refered to as the larger maximal disc on which $f$ is one-to-one.} $\Delta_f$,
i.e.\@ the maximal disc $U$,
about the origin, such that the full conjugacy
$g\circ f\circ g^{-1}(x)=\lambda x$,
holds for all $x\in U$. Estmates of linearization discs have appeard in several papers concerning the $p$-adic case  
\cite{Ben-Menahem:1988,Thiran/EtAL:1989,Pettigrew/Roberts/Vivaldi:2001,Khrennikov:2001a} later generalized in \cite{Lindahl:2004cpp}, and in \cite{Lindahl:2004,Lindahl:2009} concerning the prime characteristic case.

In this paper we consider the reaming case, namley the case in which both $K$ and the associated residue field $\Bbbk$ are of charactersitic zero. For example, $K$ could be the function field $\mathbb{C}((T))$; the field of formal Laurent series in variable $T$ over the complex numbers. Our main result can be stated in the following way.

\begin{theorem}
Let char $K=$ char $\Bbbk =0$ and $f\in K[[x]]$ be of the form (\ref{eq form of f}) with $a=\sup_{i\geq 2}|a_i|^{1/(i-1)} $. Let $\overline{\lambda}$ be the representative of $\lambda $ in $\Bbbk$ and $\Gamma (\Bbbk )$ be the set of roots of unity in $\Bbbk$. Then, the corresponding linearization disc $\Delta_f$ can be estimated as follows.
\begin{enumerate}[1.]

\item If $\overline{\lambda}\notin \Gamma (\Bbbk )$, then   $D_{1/a}(0)\subseteq \Delta_f \subseteq \overline{D}_{1/a}(0)$. If in addition $a=\max_{i\geq 2}
   |a_i|^{1/(i-1)}$ is attained as for polynomials,  or $f$ diverges on the sphere $S_{1/a}(0)$, then $\Delta_f = D_{1/a}(0)$.

 \item
  If $\overline{\lambda}\in \Gamma (\Bbbk )$ and $m$ is the smallest integer
  such that $|1-\lambda ^m|<1$. Then, \\ $ D_{\rho}(0)\subseteq\Delta_f\subseteq \overline{D}_{1/a}(0)$, where $\rho=\sqrt[m]{|1-\lambda
  ^m|}/a$. If $a=\max_{i\geq 2}
   |a_i|^{1/(i-1)}$ or $f$ diverges on the sphere $S_{1/a}(0)$, then 
  $ D_{\rho}(0)\subseteq\Delta_f\subseteq D_{1/a}(0)$.

\end{enumerate}
These estimates are maximal in the sense that there exist examples
of such $f$ which have a periodic point on the sphere
$S_{\rho}(0)$, breaking the conjugacy there.
\end{theorem}
These estimates take a simpler form than the corresponding $p$-adic case  \cite{Lindahl:2004cpp}. Second, the radius of the linearization disc is in general larger than in the 
$p$-adic case. These two facts both stem from the fact that the geometry of the roots of unity in $\mathbb{C}_p$ is more complex than in the equal characteristic case.

\section{Preliminaries}
Throughout this paper $K$ is a ultrametric field of characteristic zero (char $K= 0$), complete with respect to a nontrivial absolute value $|\cdot |$. That is, $|\cdot |$ is a multiplicative function
from $K$ to the nonnegative real numbers with $|x|=0$
precisely when $x=0$, satisfying the  
following
strong or ultrametric triangle inequality:
\begin{equation}\label{sti}
|x+y| \leq  \max[|x|,|y|],\quad\text{for all $x,y\in K$},
\end{equation}
and nontrivial in the sense that it is not
identically $1$ on $K^*$, the set of all nonzero elements
in $K$. 
One useful consequence of ultrametricity is that for any $x,y\in
K$ with $|x|\neq |y|$, the inequality (\ref{sti}) becomes an
equality. In other words, if $x,y\in K$ with $|x|<|y|$, then
$|x+y|=|y|$.

In this context it is standard to denote by $\mathcal{O}$, the ring of integers of $K$, given by
$\mathcal{O}=\{x\in K : |x|\leq 1\}$, by $\mathcal{M}$ the unique maximal ideal of
$\mathcal{O}$, given by $\mathcal{M}=\{x\in K: |x|<1\}$, and by $k$ the corresponding \emph{residue field} 
\[
\Bbbk =\mathcal{O}/\mathcal{M}.
\]
Note that if $x,y\in\mathcal{O}$ reduce to \emph{residue classes}
$\overline{x},\overline{y}\in k$, then $|x-y|$ is $1$ if $\overline{x}\neq\overline{y}$, and it is
strictly less than $1$ otherwise. 
Note also that if $K$ has positive characteristic $p$, then also char $\Bbbk =p$;  but if char $K=0$, then $\Bbbk$ could have characteristic $p$ (the $p$-adic case) or $0$. In this paper we mainly consider the latter, equal charactersitic case char $K=$ char $\Bbbk =0$. 

\begin{example}[char $K=$ char $\Bbbk =0$]\label{example laurent series char 0}
Let $F$ be a field of characteristic zero, e.g.\@ $F$ could be either $\mathbb{Q}$, $\mathbb{Q}_p$,
$\mathbb{R}$ or $\mathbb{C}$. 
Let $F((T))$be the field of
all formal Laurent series in variable $T$, with coefficients
in the field $F$.  An element $x\in K$ is of the form
\begin{equation}\label{laurent series}
x=\sum_{i\geq i_0} x_iT^i, \quad x_{i_0}\neq 0, \text{ } x_i\in
F,
\end{equation}
for some integer $i_0\in \mathbb{Z}$. 
Given $0<\epsilon <1$, we define an absolute value $|\cdot |$ on $K$ such that $|T|=\epsilon$ and
\begin{equation}\label{definition absolute value}
\left |\sum_{i\geq i_0} x_iT^i\right |=\epsilon ^{i_0}.
\end{equation}
Hence, $\pi=T$ is a uniformizer of $K$. Furthermore, $K$ is complete with respect to $|\cdot |$ and, analogously to the $p$-adic numbers, can be viewed as the completion of the field of
rational functions $F(T)$ over $F$ with respect to the absolute
value defined by (\ref{definition absolute value}). Note that in this case the residue field $\Bbbk =F$.

Note also that $i_0$ is the order of the zero (or if negative, the order of the pole) of $x$ at $T=0$.   Moreover,  $|\cdot |$ is the
trivial absolute value on $F$, the subfield of $K$ consisting of all
constant series in $K$. As for the $p$-adic  numbers, we can construct a completion $\widehat{K}$ of an algebraic closure
of $K$ with respect to an extension of $|\cdot |$. 
Then $\widehat{K}$ is a complete, algebraically closed ultrametric field, and its
residue field $\hat{\Bbbk}$ is an algebraic closure of $\Bbbk=F$. It follows that
$\hat{\Bbbk}$ has to be infinite. The value group $|\widehat{K}^*|$, that is the set of real numbers which
are actually  absolute values of non-zero elements of $\widehat{K}$, will consist
of all rational powers of $\epsilon$, rather than just integer powers of $\epsilon$ as in $|K^*|$.
In particular, the absolute value is discrete on $K$ but not on $\widehat{K}$.   
\end{example}

We use the following notation for discs.
Given an element $x\in K$ and real number $r>0$ we denote by $D_{r}(x)$ the open
disc of radius $r$ about $x$, by $\overline{D}_r(x)$ the closed
disc, and by $S_{r}(x)$ the sphere of radius $r$ about $x$.
To omit confusion, we sometimes write $D_r(x,K)$ rather than $D_{r}(x)$
to emphasize that the disc is considered as a disc in $K$.  

If $r\in|K^*|$ (that is if $r$ is actually the absolute value of some
nonzero element of $K$), we say that $D_{r}(x)$ and
$\overline{D}_r(x)$ are \emph{rational}. Note that $S_r(x)$ is
non-empty if and only if $\overline{D}_r(x)$ is rational. If
$r\notin |K^*|$, then we will call $D_{r}(x)=\overline{D}_r(x)$ an
\emph{irrational} disc. In particular, if $a\in K\subset
\mathbb{C}_p$ and $r=|a|^s$ for some rational number
$s\in\mathbb{Q}$, then $D_{r}(x)$ and $\overline{D}_r(x)$ are
rational considered as discs in the algebraic closure
$\mathbb{C}_p$. However, they may be irrational considered as
discs in $K$. Note that all discs are both open and closed as
topological sets, because of ultrametricity.
However, as we will see in Section \ref{section ultrametric
power series} below, power series distinguish between rational
open, rational closed, and irrational discs.

\subsection{Mapping properties}\label{section ultrametric power series}

Let $f$ be a power series over $K$ of the form
\begin{equation*}
f(x)=\sum_{i=0}^{\infty}a_i(x-\alpha )^i, \quad a_i\in K.
\end{equation*}
Then, $f$ converges on the open disc $D_{R_f}(\alpha )$ of radius
\begin{equation}\label{radius of convergence}
R_f = \frac{1}{\limsup |a_i| ^{1/i}},
\end{equation}
and diverges outside the closed disc $\overline{D}_{R_f}(\alpha )$
in $K$. The power series $f$ converges on the sphere
$S_{R_f}(\alpha )$ if and only if
\[
\lim_{i\to\infty}|a_i| R_f ^i=0.
\]

The basic mapping properties of ultrametric power series on discs
are given by the following generalization by Benedetto \cite{Benedetto:2003a},
of the  Weierstrass Preparation Theorem \cite{BoschGuntzerRemmert:1984,FresnelvanderPut:1981,
Koblitz:1984}.

\begin{proposition}[Lemma 2.2 \cite{Benedetto:2003a}]\label{proposition-discdegree}
Let $K$ be algebraically closed. Let
$f(x)=\sum_{i=0}^{\infty}a_i(x-\alpha)^i$ be a nonzero power
series over $K$ which converges on a rational closed disc
$U=\overline{D}_R(\alpha)$, and let $0<r\leq R$. Let
$V=\overline{D}_r(\alpha)$ and $V'=D_r(\alpha)$. Then
  \begin{eqnarray*}
    s &=& \max\{|a_i|r^i:i\geq 0\},\\
    d &=& \max\{i\geq 0:|a_i|r^i=s\},\quad and\\
    d'&=& \min\{i\geq 0:|a_i|r^i=s\}
  \end{eqnarray*}
are all attained and finite. Furthermore,
\begin{enumerate}[a.]
 \item $s\geq |f'(x_0)|\cdot r$.
 \item if $0\in f(V)$, then $f$ maps $V$ onto $\overline{D}_s(0)$
 exactly $d$-to-1 (counting multiplicity).
 \item if $0\in f(V')$, then $f$ maps $V'$ onto $D_s(0)$
 exactly $d'$-to-1 (counting multiplicity).
\end{enumerate}
\end{proposition}

We will consider the case $a_0=0$ in more detail. For our purpose,
it is then often more convenient to state Proposition
\ref{proposition-discdegree} in the following way.

\begin{proposition}\label{proposition one-to-one}
Let $K$ be algebraically closed and let
$h(x)=\sum_{i=1}^{\infty}c_i(x-\alpha )^i$ be a power series over
$K$.
\begin{enumerate}[1.]

 \item Suppose that $h$ converges on the rational closed disc
  $\overline{D}_R(\alpha)$. Let $0<r\leq R$ and suppose that
  \begin{equation}\label{ck inequality one-to-one}
   |c_i|r^i\leq |c_1|r\quad \text{ for all } i\geq 2 .
  \end{equation}
 Then, $h$ maps the open disc $D_{r}(\alpha )$ one-to-one onto
 $D_{|c_1|r}(0)$. Furthermore, if
 \[
 d = \max\{i\geq 1:|c_i|{r}^i=|c_1| r\},
 \]
 then $h$ maps the  closed disc $\overline{D}_{r}(\alpha )$ onto
 $\overline{D}_{|c_1|r}(0)$ exactly $d$-to-1 (counting
 multiplicity).

 \item Suppose that $h$ converges on the rational open disc
  $D_R(\alpha )$ (but not necessarily on the sphere $S_R(0)$).
  Let $0<r\leq R$ and suppose that
  \[
   |c_i|r^i \leq |c_1|r\quad \text{ for all } i\geq 2 .
  \]
  Then, $h$ maps $D_{r}(\alpha )$ one-to-one
  onto $D_{|c_1|r}(0)$.

\end{enumerate}

\end{proposition}

Now, suppose that $f$ has a fixed point at $x_0$ (so that $f(x_0)=x_0$) and that $|f'(x_0)|=1$. 
As a consequence of Proposition \ref{proposition-discdegree}, $f$ is not only one-to-one but a bijective  isometry on some non-empty disc about $x_0$. 
The maximal such disc is given by the
following proposition.

\begin{proposition}\label{proposition  f one-to-one}
Let $K$ be algebraically closed. Let $f\in K[[x]]$ be convergent
on some non-empty disc about $x_0\in K$. Suppose that $f(x_0)=x_0$
and $|f'(x_0)|=1$, and write
\[
f(x)=x_0+\lambda (x-x_0)+\sum_{i\geq 2}a_i(x-x_0)^i, \quad
a=\sup_{i\geq 2}|a_i|^{1/(i-1)}.
\]
Let $M$ be the largest disc, with $x_0\in M$, such that $f:M\to M$
is bijective (and hence isometric). Then $M=D_{1/a}(x_0)$ if
either $\max_{i\geq 2 }|a_i|^{1/(i-1)}$ is attained (as for polynomials) or $f$
diverges on $S_{1/a}(x_0)$. Otherwise,
$M=\overline{D}_{1/a}(x_0)$.
\end{proposition}
A proof is given in \cite{Lindahl:2004cpp}. It follows that if $f$ converges on the
sphere $S_{1/a}(x_0)$ but fails to be one-to-one there, then there
is a point $x\in S_{1/a}(x_0)$ such that $f(x)=x_0=f(x_0)$. This
is always the case when $f$ is a polynomial.

That $f$ may diverge on $S_{1/a}(x_0)$ follows since, for example, the
power series $f(x)=\lambda x + \sum_{i=2}^{\infty}(a_2)^{i-1}x^i$
converges if and only if $|x|<1/|a_2|=1/a$.

Furthermore, for every $x\in M$, $|f(x)-x_0|=|x-x_0|$ and hence
all spheres in $M$ are invariant under $f$.

\begin{remark}\label{remark Lemma bijective}
Recall that the discs $D_{1/a}(0)$ and $\overline{D}_{1/a}(0)$ are
rational if and only if $a=\sup_{i\geq 2}|a_i|^{1/(i-1)}\in |K|$.
If the maximum $a=\max_{i\geq 2}|a_i|^{1/(i-1)}$ exists, and $K$
is algebraically closed,  then $a\in|K|$. This is always the case
if $f$ is a polynomial. If $f$ is not a polynomial and the maximum
fails to exist we may have $\sup_{i\geq 2} |a_i|^{1/(i-1)}\notin
|K|$. Let $K=\mathbb{C}_p$. Let $\beta $ be an irrational number
and let $p_n/q_n$ be the $n$-th convergent of the continued
fraction expansion of $\beta$. Let the sequence $\{a_i\in
\mathbb{Q}_p\}_{i\geq 2}$ satisfy
\[
|a_i|= \left \{
\begin{array}{ll}
p^{p_n}, & \textrm{if \quad $i-1=q_n$ and $p_n/q_n<\beta $},\\
0, & \textrm{otherwise}.
\end{array}\right.
\]
Then,
\[
\sup_{i\geq 2} |a_i|^{1/(i-1)}=p^{\beta}\notin
|K|=\{p^r:r\in\mathbb{Q} \}\cup \{0\}.
\]
\end{remark}

For more information on ultrametric power series the reader
can consult
\cite{Schikhof:1984}.
From a dynamical point of view,  the paper \cite{Benedetto:2003a}
contains many useful results on ultrametric analogues of
complex analytic mapping theorems relevant for dynamics.

\subsection{The linearization disc}

The results above have some important implications for linearization 
discs. We use the following definition of a linearization disc. Let $K$
be a complete ultrametric field of charactersitic zero. Suppose that $f\in K[[x]]$
has an indifferent fixed point $x_0\in K$, with multiplier
$\lambda =f'(x_0)$, not a root of unity. By
\cite{Herman/Yoccoz:1981}, there is a unique formal power series
solution $g$, with $g(x_0)=0$ and $g'(x_0)=1$, to the following
form of the Schr\"{o}der functional equation
\[
g\circ f(x)=\lambda g(x).
\]
If the formal solution $g$ converges on some non-empty disc about
$x_0$, then the corresponding \emph{linearization disc} of $f$ about
$x_0$, denoted by $\Delta _f(x_0)$, is defined as the largest disc
$U\subset K$, with $x_0\in U$, such that the Schr\"{o}der
functional equation holds for all $x\in U$, and $g$ converges and
is one-to-one on $U$. We will often refer to $g$ as the
\emph{conjugacy function}.

This notion of a linearization disc is well-defined since, by proposition \ref{proposition f one-to-one},
there always exist a largest disc on which $g$ is one-to-one (provided that $g$
is convergent). Recall that by the ultrametric Siegel
theorem by Herman and Yoccoz \cite{Herman/Yoccoz:1981}, the formal solution $g$ always
converges if char $K=0$.

As a consequence of the results stated in the section above, both $f$ and the
conjugacy $g$ turn out to be one-to-one and isometric on a
ultrametric linearization disc.

\begin{proposition}\label{proposition linearization disc isometry}
Let $K$ be algebraically closed. Suppose that $f\in K[[x]]$ has a
linearization disc $\Delta_f(x_0)$ about $x_0\in K$. Let $g$, with
$g(x_0)=0$ and $g'(x_0)=1$, be the corresponding conjugacy
function. Then, both $g:\Delta_f(x_0)\to g(\Delta_f(x_0))$ and
$f:\Delta_f(x_0)\to \Delta_f(x_0)$ are bijective and isometric. In
particular,
if $x_0=0$, then $g(\Delta_f(x_0))=\Delta_f(x_0)$. Furthermore,
$\Delta_f(x_0)\subseteq M\subseteq \overline{D}_{1/a}(x_0)$, where
$M$ and $a$  are defined as in Lemma \ref{proposition  f one-to-one}.
\end{proposition}

Hence, the radius of a linearization disc $\Delta_f(x_0)$ is
equal to to that of $g(\Delta_f(x_0))$. In particular, the radius
of a linearization disc is independent of the location of the fixed point
$x_0$. Therefore, we shall, without loss of generality, henceforth
assume that $x_0=0$.

\begin{remark}
All the results in this and the previous section, except for Proposition
\ref{proposition-discdegree}, hold also in the case that $K$ is
not algebraically closed, with the modification that the mappings
are are one-to-one but not necessarily surjective. 
\end{remark}

\section{Proof of the main theorem}
From now on char $K=$ char $\Bbbk =0$.
As noted in the previous section, we may, without loss of
generality,  assume that $f$ has its fixed point at the origin,
and that $f\in \mathcal{F}_{\lambda,a}$, as defined below. Let
$\lambda\in K$ be such that
\begin{equation}
|\lambda|=1, \quad\textrm{but } \lambda ^n\neq 1,\quad  \forall
n\geq 1,
\end{equation}
and let $a$ be a real number. We shall associate with the pair
$(\lambda,a)$ a family $\mathcal{F}_{\lambda,a}$ of power series
defined by
\begin{equation}
\mathcal{F}_{\lambda,a}:=\left\{\lambda x
+\sum a_ix^i\in\mathbb{C}_p[[x]]:a=\sup_{i\geq 2}|a_i|^{1/(i-1)}\right\}.
\end{equation}
It follows that each $f\in \mathcal{F}_{\lambda,a}$ is convergent
on $D_{1/a}(0)$, and by Proposition \ref{proposition  f one-to-one}
$f:D_{1/a}(0)\to D_{1/a}(0)$ is bijective and isometric.

As $K$ is of characteristic zero, we may, by the
ultrametric Siegel theorem \cite{Herman/Yoccoz:1981},
associate with $f$ a unique convergent power series solution $g$
to the Scr\"{o}der functional equation, of the form
\[
g(x)=x + \sum_{k\geq 2}b_kx^k,
\]
and a corresponding linearization disc about the origin
\[
\Delta_f:=\Delta _f(0).
\]
Recall that by Proposition \ref{proposition linearization disc isometry}, since
$x_0=0$, the linearization disc $\Delta_f$ is the largest disc $U\subset
K$ about the origin such that the full conjugacy
$g\circ f \circ g^{-1}(x)=\lambda x$ holds for all $x\in U$.

Given $f\in \mathcal{F}_{\lambda,a}$, Proposition 
\ref{proposition linearization disc isometry} yields the following concerning $\Delta_f$.

\begin{lemma}\label{lemma upper bound Siegel and isometry}
Let $f\in \mathcal{F}_{\lambda,a}$. Then $f$ has a linearization disc
$\Delta_f$ about the origin in $K$. Let $g$, with
$g(0)=0$ and $g'(0)=1$, be the corresponding conjugacy function.
Then, the following two statements hold:
\begin{enumerate}[1)]

   \item Both $g:\Delta_f\to \Delta_f$ and $f:\Delta_f\to \Delta_f$ are bijective
   and isometric.

   \item $\Delta_f\subseteq \overline{D}_{1/a}(0)$. If $a=\max_{i\geq 2}
   |a_i|^{1/(i-1)}$ or $f$ diverges on the sphere $S_{1/a}(0)$, then $\Delta_f\subseteq D_{1/a}(0)$.

\end{enumerate}
\end{lemma}

Also note the following lemma concerning estimates of the coefficients of the conjugacy.

\begin{lemma}\label{lemma bk estimate indiff}
Let $f\in \mathcal{F}_{\lambda,a}$. Then, the coefficients of the
conjugacy function $g$ satisfy
\begin{equation}\label{b_k estimate by prod 1-lambda n}
|b_k|\leq \left( \prod_{n=1}^{k-1}|1-\lambda ^n| \right
)^{-1}a^{k-1},
\end{equation}
for all $k\geq 2$.
\end{lemma}

\begin{proof} The coefficients of the conjugacy $g$ must
satisfy the recurrence relation
\begin{equation}\label{bk-equation}
b_k=\frac{1}{\lambda (1-\lambda^{k-1})}\sum_{l=1}^{k-1}b_l
(\sum\frac{l!}{\alpha_1!\cdot ...\cdot
\alpha_k!}a_1^{\alpha_1}\cdot ...\cdot a_k^{\alpha_k})
\end{equation}
where $\alpha _1,\alpha _2,\dots,\alpha _k$ are nonnegative
integer solutions of
\begin{equation}\label{index-equations}
   \left\{\begin{array}{ll}
            \alpha_1+...+\alpha_k=l,\\
            \alpha_1+2\alpha_2...+k\alpha_k=k,\\
            1\leq l\leq k-1.
        \end{array}
   \right.
\end{equation}
Note that the factorial factors
$l!/\alpha_1!\cdot\dots\cdot\alpha_k!$ are always integers and
thus of modulus less than or equal to $1$. Also recall that
$|a_i|\leq a^{i-1}$.  It follows that
\[
|b_k|\leq \left( \prod_{n=1}^{k-1}|1-\lambda ^n| \right )^{-1}
a^{\alpha},
\]
for some integer $\alpha$. In view of equation
(\ref{index-equations}) we have
\[
\sum_{i=2}^{k}(i-1)\alpha_i=k-l.
\]
Consequently, since $|a_i|\leq a^{i-1}$, we obtain
\begin{equation}\label{estimate prod a_i}
\prod_{i=2}^k|a_i|^{\alpha_i}\leq
\prod_{i=2}^{k}a^{(i-1)\alpha_i}=a^{k-l}.
\end{equation}
Now we use induction over $k$. By definition $b_1=1$ and,
according to the recursion formula (\ref{bk-equation}), $|b_2|\leq
|1-\lambda |^{-1}|a_2|\leq |1-\lambda |^{-1}|a|$. Suppose that
\[
|b_l|\leq \left( \prod_{n=1}^{l-1}|1-\lambda ^n| \right
)^{-1}a^{l-1}
\]
for all $l<k$. Then
\[
|b_k|\leq \left( \prod_{n=1}^{k-1}|1-\lambda ^n| \right
)^{-1}a^{l-1}\max\left\{ \prod_{i=2}^k|a_i|^{\alpha_i}\right \},
\]
and the lemma follows by the estimate (\ref{estimate prod a_i}).
\end{proof}

In the following we show how to calculate the distance
$|1-\lambda ^n |$ for an arbitrary integer $n\geq 1$. Applying
Proposition \ref{proposition one-to-one} to the estimate in the
above lemma we can then estimate the disc on which the conjugacy
function $g$ is one-to-one.

Let $\lambda\in S_1(0)$, be an element in the unit sphere. We are
concerned with calculating the distance
\[
|1-\lambda ^n|, \quad\text{for $n=1,2,\dots$}.
\]
Recall that if $x,y\in\overline D_1(0)$, then  $|x-y|<1$ if and only
if the reductions $\overline{x},\overline{y}$ belong to the same
residue class. Consequently,
\begin{equation}\label{distance and residue}
|1-\lambda ^n|<1\quad\iff\quad\overline{\lambda}^n -
1=0\quad\text{in $\Bbbk$}.
\end{equation}
Hence, the behavior of $1-\lambda ^n$ falls into one of two
categories, depending on whether the reduction of $\lambda $ is a
root of unity or not. For convenience, denote by $\Gamma (\Bbbk)$ the set of roots of unity in 
$\Bbbk$. More precisely, we have the following lemma.

\begin{lemma}\label{lemma distance char 0,0}
Let char $K=$ char $\Bbbk=0$. Suppose $\lambda\in S_1(0)$.
Then
\begin{enumerate}
  \item $\overline{\lambda}\notin\Gamma(\Bbbk)$ $\iff$ $|1-\lambda ^n|=1$
  for all integers $n\geq 1$.
  \item If $\overline{\lambda}\in\Gamma(\Bbbk)$ then there is a smallest  integer $m\geq
  1$ such that $|1-\lambda ^m|<1$. Moreover,
\begin{equation}\label{distance char 0}
  \left |1- \lambda ^n \right |=
  \left\{\begin{array}{ll}
         1, & \textrm{ if \quad $m\nmid n$,}\\
        |1-\lambda ^m|, & \textrm{ if \quad $m\mid n$.}
  \end{array}\right.
\end{equation}

\end{enumerate}

\end{lemma}
\begin{proof} First note that since the characteristic of
the residue field is zero, $ | \binom{l}{k} |=1$ for all binomials
$\binom{l}{k}$. By ultrametricity (\ref{sti}) we then have
\[
|((\lambda ^m-1)+1)^l-1|=|\sum_{k=1}^{l}\binom{l}{k}(\lambda ^m
-1)^k|=|\lambda ^m-1|,
\]
as required.  
\end{proof}

\vspace{1.5ex} \noindent Note that category 2 in Lemma \ref{lemma
distance char 0,0} is always non-empty since $1,-1 \in
\mathbb{Q}\subseteq \Bbbk $.

With these results at hand, we are now in a position to prove our main result. 

\begin{theorem}
Let char $K = $ char $\Bbbk=0$ and $f\in\mathcal{F}_{\lambda,a}$. Then, the corresponding linearization disc $\Delta_f$ can be estimated as follows.
\begin{enumerate}[1.]

\item If $\overline{\lambda}\notin \Gamma (\Bbbk )$, then   $D_{1/a}(0)\subseteq \Delta_f \subseteq \overline{D}_{1/a}(0)$. If, in addition, $a=\max_{i\geq 2}
   |a_i|^{1/(i-1)}$ as for polynomials,  or $f$ diverges on the sphere $S_{1/a}(0)$, then $\Delta_f = D_{1/a}(0)$.

 \item
  If $\overline{\lambda}\in \Gamma (\Bbbk )$ and $m$ is the smallest integer
  such that $|1-\lambda ^m|<1$. Then, \\ $ D_{\rho}(0)\subseteq\Delta_f\subseteq \overline{D}_{1/a}(0)$, where $\rho=\sqrt[m]{|1-\lambda
  ^m|}/a$. If $a=\max_{i\geq 2}
   |a_i|^{1/(i-1)}$ or $f$ diverges on the sphere $S_{1/a}(0)$, then 
  $ D_{\rho}(0)\subseteq\Delta_f\subseteq D_{1/a}(0)$.

\end{enumerate}
These estimates are maximal in the sense that there exist examples
of such $f$ which have a periodic point on the sphere
$S_{\rho}(0)$, breaking the conjugacy there.
\end{theorem}

\begin{proof}
The first statement follows immediately from Lemma \ref{lemma upper bound Siegel and isometry}, Lemma \ref{lemma bk estimate indiff}, and Lemma \ref{lemma distance char 0,0}.
To see that $f$ may have a periodic point on the boundary breaking the conjugacy there we consider the polynomial $f(x)=\lambda x +
a_nx^n$ for some integer $n\geq 2$. Note that we can choose
$a=|a_n|^{1/(n-1)}$ in this case. Moreover $\hat{x}=[(1-\lambda
)/a_n]^{1/(n-1)}$ is fixed under $f$. But, since $|1-\lambda |=1$,
we have $|\hat{x}|=1/a$ so that $\hat{x}$ sits on the sphere
$S_{1/a}(0)$.

Now we consider the second statement. In this case, Lemma
\ref{lemma distance char 0,0} implies that
\[
|b_k|\leq \frac{1}{|1-\lambda ^m|^{\lfloor k-1/m \rfloor}}a^{k-1}.
\]
It follows that the conjugacy $g$ converges on the open disk of
radius
\[
\left ( \limsup |b_k|^{1/k}\right )^{-1}\geq |1-\lambda
^m|^{\frac{1}{m}}a^{-1}=\rho.
\]
Moreover,
\[
|b_k|\leq |1-\lambda ^m|^{-\frac{k-1}{m}}a^{k-1}=\rho ^{-(k-1)},
\]
Consequently,
\[
|b_k| \rho ^k \leq \rho =|b_1|\rho,
\]
In view of Proposition \ref{proposition one-to-one},
$g:D_{\rho}(0)\to D_{\rho}(0)$ is bijective. Recall that by Lemma
\ref{proposition  f one-to-one} $f:D_{1/a}(0)\to D_{1/a}(0)$ is a
bijection. Moreover, $1/a>\rho$.  Consequently, the linearization disc 
$\Delta_f\supseteq D_{\rho}(0)$ as required.
This estimate of $\Delta_{f}$ is maximal in the
sense that all $f$ of the form $f(x)=\lambda x+ a_2x^2$ have a
periodic fixed point on the sphere $S_{\rho}(0)$. 
\end{proof}

\section*{Acknowledgements}

I would like to thank
Prof. Andrei Yu. Khrennikov for fruitful discussions and introducing me to the
theory of ultrametric dynamical systems and conjugate maps.

\end{document}